\documentclass[12pt]{article}

\usepackage{latexsym, amssymb, amsmath, amscd, amsfonts, epsfig, graphicx, colordvi,verbatim,ifpdf}
\usepackage{amsfonts, amsmath, amssymb}
\usepackage{amssymb,amsfonts,amsmath,latexsym,epsfig,cite, psfrag,eepic,color}
\usepackage{amscd,graphics}
\usepackage{latexsym, amssymb,  amsmath,amscd, amsfonts, epsfig, graphicx, colordvi,amsthm}
\usepackage{float}

\usepackage{graphicx}
\usepackage{color}
\usepackage{ifpdf}
\usepackage{fancybox}
\usepackage{makecell}
\usepackage{multicol}
\usepackage{multirow}

\newtheorem{thm}{Theorem}[section]

\newtheorem{rem}[thm]{\it Remarks}
\newtheorem{defi}[thm]{Definition}
\newtheorem{lem}[thm]{Lemma}

\newtheorem{exa}[thm]{Example}

\def\qed{\nopagebreak\hfill{\rule{4pt}{7pt}}
\medbreak}

\def\qed{\nopagebreak\hfill{\rule{4pt}{7pt}}
\medbreak}

\numberwithin{equation}{section}

\def\qed{\nopagebreak\hfill{\rule{4pt}{7pt}}
\medbreak}

\newlength{\boxedparwidth}
  {\begin{center} \begin{tabular}{|@{\hspace{.315in}}c@{\hspace{.15in}}|}
                  \hline \\ \begin{minipage}[t]{\boxedparwidth}
                  \setlength{\parindent}{.25in}}%
  {\end{minipage} \\ \\ \hline \end{tabular} \end{center}}

\begin{document}
\begin{center}

{ \large\bf Combinatorial proofs of the Ramanujan type congruences modulo 3}
\end{center}

\vskip 5mm

\begin{center}
{  Robert X. J. Hao}\\
    Department of Mathematics and Physics, \\Nanjing Institute of Technology,
    Nanjing 211167, P.~R.~China\\[6pt]

    haoxj@njit.edu.cn
\end{center}

\vskip 6mm \noindent {\bf Abstract.}
The partition statistic $V_R$-rank is introduced to give combinatorial proofs of the Ramanujan type congruences mod 3 for certain classes of partition functions.

\vskip 5mm
\noindent {\bf Keywords}: partition statistic; bijection; combinatorial proof; Ramanujan type congruence;

\noindent {\bf AMS Classifications}: 05A17, 11P83

\section{Introduction}

A partition $\lambda$ of a positive integer $n$
is a weakly-decreasing sequence of positive integers $\lambda_1\geq\lambda_2\geq\ldots\geq\lambda_l$ such that
$|\lambda|=\sum_{i=1}^{l}\lambda_i=n$. Here we adopt the convention that the partition of 0 is the empty set.
Let $p(n)$ denote the number of
partitions of $n$.
In 1919, Ramanujan \cite{Ramanujan-1919} found the following three congruences on $p(n)$
\begin{align}
\label{eqnr1}p(5n+4)&\equiv0 \pmod{5}, \\[5pt]
\label{eqnr2}p(7n+5)&\equiv0 \pmod{7}, \\[5pt]
\label{eqnr3}p(11n+6)&\equiv0\pmod{11}.
\end{align}
In 1988, Andrews and Garvan \cite{Andrews-Garvan-1988,Garvan-1988} defined the partition
statistic ``crank'', and provided combinatorial interpretations for congruences \eqref{eqnr1}--\eqref{eqnr3}.

Later, Garvan, Kim, and Stanton \cite{Garvan-Kim-Stanton-1990} found different cranks, which also explained all three congruences. It should be mentioned that they established explicit bijections between various equinumerous classes, and
in particular, provided a combinatorial proof of \eqref{eqnr1}.

By imposing certain restrictions on the parts of the partitions, one can obtain variants of restricted partitions.
For instance, an overpartition is a partition for which the first occurrence of a part may be overlined and a partition with designated summands is a partition formed by tagging exactly one part among parts with equal size.
Properties of these partition functions have been investigated quite intensively
in recent studies, see, for example
\cite{Andrews-Lewis-Lovejoy-2002,Hirschhorn-Sellers-2005,Chen-Ji-Jin-Shen-2013,
Fu-Tang-2017,Wang-2017,Wang-Yee-2020
}.
Suppose that $t$ is a positive integer and $S$ is a set of non-negative integers. Let us recall variants of significant restricted partition sets given as
\begin{align*}
\mathcal{P}_{t_S}=&\{\lambda|s\in S, \lambda \ \textrm{is a partition into parts $\equiv s \pmod{t}$}\},\\[4pt]
\mathcal{D}_{t_S}=&\{\lambda|s\in S, \lambda \ \textrm{is a partition into distinct parts $\equiv s \pmod{t}$}\},\\[4pt]
\overline{\mathcal{P}}=&\{\lambda| \lambda\  \textrm{is an overpartition}\},\\[4pt]
\mathcal{S}=&\{\lambda| \lambda \ \textrm{is a staircase partition}\},\\[4pt]
\mathcal{S}_{t-odd}=&\{\lambda| \lambda \ \textrm{is a staircase partition into odd parts, with the part 1 can be overlined}\},\\[4pt]
\mathcal{PD}=&\{\lambda| \lambda \ \textrm{is a partition with designated summands}\},\\[4pt]
\mathcal{POD}=&\{\lambda| \lambda \ \textrm{is a partition with odd parts distinct}\},\\[4pt]
\mathcal{A}=&\left\{\lambda| \lambda\  \textrm{is a 2-color partition colors $r$ and $b$}\right.\\[4pt] &\left.\textrm{such that the color $b$ appears only in even parts}\right\}.
\end{align*}
For example, it is easy to see that $\mathcal{P}_{1_{\{0\}}}$ denotes the set of
ordinary partitions, and $\mathcal{D}_{2_{\{1\}}}$ denotes the set of partitions into distinct odd parts.
Correspondingly, for a given set of partitions denoted by a script capital letter (or letters), e.g., $\mathcal{PD}$, we denote the script capital letters $\mathcal{PD}^n$ the set of partitions $\lambda\in \mathcal{PD}$ with $|\lambda|=n$, and denote the minuscules $pd(n)$ the number of elements in $\mathcal{PD}^n$.

For convenience, one may consider a special kind of restricted partition as $k$-colored partition, which are $k$ ordered restricted or unrestricted partitions $\lambda^{(1)},\lambda^{(2)},\dots,\lambda^{(k)}$.
It is not hard to check that $\lambda \in \mathcal{A}$ can be seen
as a 2-colored partition with $\lambda^{(1)} \in \mathcal{P}_{1_{\{0\}}} ,\lambda^{(2)} \in \mathcal{P}_{2_{\{0\}}}$. A special but important case is that $\lambda^{(1)},\lambda^{(2)},\dots,\lambda^{(k)}$ subject to the same restrictions, for example
\begin{align*}
{\mathcal{P}_{-k}}=&\left\{\overrightarrow{\lambda}
=(\lambda^{(1)},\lambda^{(2)},\dots,\lambda^{(k)})|\lambda^{(i)} \in \mathcal{P}_{1_{\{0\}}}\  \textrm{for} \ 1\leq i\leq k\right\},\\[4pt]
\overline{\mathcal{P}}_{-k}=&\left\{\overrightarrow{\lambda}
=(\lambda^{(1)},\lambda^{(2)},\dots,\lambda^{(k)})|\lambda^{(i)} \in \overline{\mathcal{P}} \ \textrm{for}\  1\leq i\leq k\right\},\\[4pt]
\mathcal{POD}_{-k}=&\left\{\overrightarrow{\lambda}
=(\lambda^{(1)},\lambda^{(2)},\dots,\lambda^{(k)})| \lambda^{(i)} \in \mathcal{POD} \ \textrm{for}\  1\leq i\leq k\right\}.
\end{align*}
By convention, we call a 2-colored partition bipartition.

In the studies of restricted partitions, some interesting Ramanujan type congruences were established, see, for example, Andrews et al. \cite{Andrews-Lewis-Lovejoy-2002}, Chan \cite{Chan-2010}, and Chen et al. \cite{Chen-Lin-2011,Chen-Lin-2012}.
\begin{align}
\label{pd-3}pd(3n+2)&\equiv0 \pmod3,\\[3pt]
\label{a-3}a(3n+2)&\equiv0 \pmod3,\\[3pt]
\label{pod-3}pod_{-2}(3n+2)&\equiv0 \pmod3,\\[3pt]
\label{o-3}\overline{p}_{-2}(3n+2)&\equiv0 \pmod3.
\end{align}
Chen, Ji, Jin, Shen \cite{Chen-Ji-Jin-Shen-2013}, Kim \cite{Kim-2010}, and Chen, Lin \cite{Chen-Lin-2011,Chen-Lin-2012} provided combinatorial interpretations of the above Ramanujan type congruences by introducing partition statistics respectively.

In this paper, we aim to give combinatorial proofs for the above congruence properties by organizing the set of restricted partitions into orbits such that each orbit consists of 3 distinct members, and each element of the orbit has a distinct partition statistic $\pmod{3}$.
To do so, we consider a special kind of $k$-colored partition and
define the partition statistic $V_{R}$-rank.
We denote
a set of restricted partitions by the script capital letter $\mathcal{R}$, and for a partition $\lambda$, denote $\ell(\lambda)$ the number of its parts.
\begin{defi}
For $k\geq3$, let
\begin{align}\label{vtk}
\nonumber V_{t,k}=&\left\{\overrightarrow{\lambda}=(\lambda^{(1)},\lambda^{(2)},\lambda^{(3)},\cdots,\lambda^{(k)})|
\lambda^{(1)},\lambda^{(2)},\lambda^{(3)} \in \mathcal{P}_{t_{\{0\}}}, \ \textrm{for}\  i>3, \right. \\[3pt]
&\left.\lambda^{(i)} \ \textrm{is an ordinary partition or a restricted partition}
\right\}.
\end{align}
\end{defi}
\begin{defi}\label{VR-rank}
For a given set of restricted partitions $\mathcal{R}$ and $\lambda \in \mathcal{R}$,
if there exist a set $V_{t,k}$ and
a bijection $\Lambda:\mathcal{R}^n \rightarrow V_{t,k}^n$ such that $ \Lambda(\lambda)=(\lambda^{(1)},\lambda^{(2)},\lambda^{(3)},\cdots,\lambda^{(k)})$,
the $V_{R}$-rank of $\lambda$, denoted by $r_{V}(\lambda)$, is defined as
\begin{align*}
r_{V}(\lambda)=\ell(\lambda^{(1)})-\ell(\lambda^{(2)}).
\end{align*}
For convenience, we define the
$V_{R}$-rank of $\overrightarrow{\lambda}\in V_{t,k}$ as
\begin{align*}
r_{V}(\overrightarrow{\lambda})=\ell(\lambda^{(1)})-\ell(\lambda^{(2)}).
\end{align*}
\end{defi}
The following theorem holds.
\begin{thm}\label{HLW-ZZ}
For a given set of restricted partitions $\mathcal{R}$,
if there exist a set $V_{t,k}$ and
a bijection $\Lambda:\mathcal{R}^n \rightarrow V_{t,k}^n$ such that for any $\overrightarrow{\lambda}\in V_{t,k}$,
\begin{equation}\label{206}
\sum_{i=4}^k|\lambda^{(i)}|\not\equiv j\pmod{3},
\end{equation}
then the set $\mathcal{R}^{3n+j}$
can be organized into orbits, with each orbit consists of 3 distinct members and each element of the orbit has a distinct $V_{R}$-rank $\pmod{3}$, where the $V_{R}$-rank is defined by Definition \ref{VR-rank}.
\end{thm}

We establish
bijections to give combinatorial proofs of \eqref{pd-3}--\eqref{o-3} based on Theorem \ref{HLW-ZZ}.

\section{Preliminaries}
In this section, we prove Theorem \ref{HLW-ZZ}.

\emph{Proof of Theorem \ref{HLW-ZZ}.}
Let $\mathcal{R}$ fulfill the conditions of Theorem \ref{HLW-ZZ}. Suppose $\lambda \in \mathcal{R}$ and $\Lambda(\lambda)=(\lambda^{(1)},\lambda^{(2)},\lambda^{(3)},\cdots,\lambda^{(k)})$.
By \eqref{206}, we have that
$\lambda \in \mathcal{R}^{3n+j}$ impiles
\begin{equation}\label{333}
|\lambda^{(1)}|+|\lambda^{(2)}|+|\lambda^{(3)}| \not\equiv 0 \pmod{3}.
\end{equation}
If $t\equiv 0 \pmod{3}$, the theorem follows immediately by the fact that
\[
|\lambda^{(1)}|+|\lambda^{(2)}|+|\lambda^{(3)}| \equiv 0 \pmod{3}.
\]
Denote $\ell_i(\lambda^{(m)})\ \ (i=\pm1,m=1,2,3)$ the number of parts $\equiv i \pmod{3}$ in $\lambda^{(m)}$.
If $t\not\equiv 0 \pmod{3}$, one can easily check
\eqref{333} equates with
\[
\ell_1(\lambda^{(1)})+\ell_1(\lambda^{(2)})+\ell_1(\lambda^{(3)})
-[ \ell_{-1}(\lambda^{(1)})+\ell_{-1}(\lambda^{(2)})+\ell_{-1}(\lambda^{(3)})] \not\equiv 0 \pmod{3}.
\]
Hence
\[
\ell_1(\lambda^{(1)})+\ell_1(\lambda^{(2)})+\ell_1(\lambda^{(3)})
\not\equiv \ell_{-1}(\lambda^{(1)})+\ell_{-1}(\lambda^{(2)})+\ell_{-1}(\lambda^{(3)}) \pmod{3}.
\]
Thus, for any $\lambda \in \mathcal{R}^{3n+j}$, there are two cases.
\begin{itemize}
\item[Case 1.] $\ell_1(\lambda^{(1)})+\ell_1(\lambda^{(2)})+\ell_1(\lambda^{(3)}) \not\equiv 0 \pmod{3};$
\item[Case 2.] $\ell_1(\lambda^{(1)})+\ell_1(\lambda^{(2)})+\ell_1(\lambda^{(3)}) \equiv 0 \pmod{3},\ell_{-1}(\lambda^{(1)})+\ell_{-1}(\lambda^{(2)})+\ell_{-1}(\lambda^{(3)}) \not\equiv 0 \pmod{3}.$
\end{itemize}
We only consider the first case,
and the second case can be justified in the same manner. For the first case,
we claim
\begin{align}
\label{lun1}\ell_1(\lambda^{(1)})-\ell_1(\lambda^{(2)})&\not\equiv \ell_1(\lambda^{(3)})-\ell_1(\lambda^{(1)}) \pmod{3},\\[4pt]
\label{lun2}\ell_1(\lambda^{(3)})-\ell_1(\lambda^{(1)})&\not\equiv \ell_1(\lambda^{(2)})-\ell_1(\lambda^{(3)})\pmod{3},\\[4pt]
\label{lun3}\ell_1(\lambda^{(2)})-\ell_1(\lambda^{(3)})&\not\equiv \ell_1(\lambda^{(1)})-\ell_1(\lambda^{(2)})\pmod{3}.
\end{align}
Take \eqref{lun1} as an example.
By
\[
\ell_1(\lambda^{(1)})+\ell_1(\lambda^{(2)})+\ell_1(\lambda^{(3)}) \not\equiv 0 \pmod{3},
\]
we have
\[
\ell_1(\lambda^{(1)})+\ell_1(\lambda^{(2)})+\ell_1(\lambda^{(3)}) -3\ell_1(\lambda^{(1)}) \not\equiv 0 \pmod{3},
\]
which implies \eqref{lun1} immediately.
The proofs of \eqref{lun2} and \eqref{lun3} are similar to that of
\eqref{lun1}, and hence they are omitted.
Invoking \eqref{lun1}-\eqref{lun3}, we see that
each element of
\[
\ell_1(\lambda^{(1)})-\ell_1(\lambda^{(2)}),\ \
\ell_1(\lambda^{(3)})-\ell_1(\lambda^{(1)}),\ \
\ell_1(\lambda^{(2)})-\ell_1(\lambda^{(3)}),
\]
has a distinct residue mod 3.
Based on this observation, we conclude that for a given $\overrightarrow{\lambda}=(\lambda^{(1)},\lambda^{(2)},\lambda^{(3)},\cdots,\lambda^{(k)})
\in V_{t,k}^{3n+j}$,
if we fix all the parts
$\not\equiv1 \pmod{3}$ in $\lambda^{(1)},\lambda^{(2)},\lambda^{(3)}$, and alternate the parts $\equiv1 \pmod{3}$ in them, we can establish $\overrightarrow{\lambda},\overrightarrow{\lambda}_{\widehat{}},
\overrightarrow{\lambda}_{\widehat{\widehat{}}}\in V_{t,k}^{3n+j}$ such that each of them has a distinct $V_{R}$-rank $\pmod{3}$.
\begin{exa}\label{example-1}
Let
\begin{align}\label{V14-o} V_{1,4}=&\left\{(\lambda^{(1)},\lambda^{(2)},\lambda^{(3)},\lambda^{(4)})|
\lambda^{(1)},\lambda^{(2)},\lambda^{(3)} \in \mathcal{P}_{1_{\{0\}}}, \lambda^{(4)} \in \mathcal{S}
\right\},
\end{align}
and $\overrightarrow{\lambda}\in V_{1,4}^{83}$ as given by
\[
\overrightarrow{\lambda}=(9+8+7+7+5+4,5+2+1,10+6+4+4+3+2,3+2+1).
\]
Since
\[
\ell_1(\lambda^{(1)})+\ell_1(\lambda^{(2)})+\ell_1(\lambda^{(3)})=7 \not\equiv 0 \pmod{3},
\]
alternating the parts $\equiv1 \pmod{3}$ in $\lambda^{(1)}$, $\lambda^{(2)}$, $\lambda^{(3)}$, we establish the following
three partitions belong to $V_{1,4}^{83}$ as defined in \eqref{V14-o}.
\begin{align*}
\overrightarrow{\lambda}=&(9+8+7+7+5+4,5+2+1,10+6+4+4+3+2,3+2+1),\\[4pt]
\overrightarrow{\lambda}_{\widehat{}}=&(10+9+8+5+4+4,7+7+5+4+2,6+3+2+1,3+2+1),\\[4pt]
\overrightarrow{\lambda}_{\widehat{\widehat{}}}=&(9+8+5+1,10+5+4+4+2,7+7+6+4+3+2,3+2+1).
\end{align*}
One can see
\begin{align*}
r_{V}(\overrightarrow{\lambda})=&6-3=3,\\[4pt]
r_{V}(\overrightarrow{\lambda}_{\widehat{}})=&6-5=1,\\[4pt]
r_{V}(\overrightarrow{\lambda}_{\widehat{\widehat{}}})=&4-5=-1.
\end{align*}
\end{exa}
To be more specific, for $i=\pm 1, m=1,2,3$, denote $\lambda_{3_{\{i\}}}^{(m)}$ the partition consisting of all the parts
$\equiv i \pmod{3}$ of $\lambda^{(m)}$, and $\lambda_{3_{\overline{\{i\}}}}^{(m)}$ the partition consisting of all the parts
$\not\equiv i \pmod{3}$ of $\lambda^{(m)}$ respectively.
For the first case, we give a new representation of $\overrightarrow{\lambda}=(\lambda^{(1)},\lambda^{(2)},\lambda^{(3)},\cdots,\lambda^{(k)})$
as follows
\begin{equation*}
\overrightarrow{\lambda}=(\overrightarrow{\lambda_1},
\overrightarrow{\lambda_{\overline{1}}};
\lambda^{(4)},\cdots,\lambda^{(k)}),
\end{equation*}
where
\begin{equation*}
\overrightarrow{\lambda_1}=
(\lambda_{3_{\{1\}}}^{(1)},\lambda_{3_{\{1\}}}^{(2)},
\lambda_{3_{\{1\}}}^{(3)}),
\end{equation*}
and
\begin{equation*}
\overrightarrow{\lambda_{\overline{1}}}=
(\lambda_{3_{\overline{\{1\}}}}^{(1)},
\lambda_{3_{\overline{\{1\}}}}^{(2)},
\lambda_{3_{\overline{\{1\}}}}^{(3)}).
\end{equation*}
It can be checked that under the following
cyclic permutation
\[
\widehat{C}_1(\overrightarrow{\lambda_1})=
(\lambda_{3_{\{1\}}}^{(3)},
\lambda_{3_{\{1\}}}^{(1)},
\lambda_{3_{\{1\}}}^{(2)})
\]
while $r_{V_3}(\overrightarrow{\lambda})$ increases or decreases by 1 $\pmod{3}$ under the map
\begin{equation*}
\widehat{O}_1(\overrightarrow{\lambda})
=(\widehat{C}(\overrightarrow{\lambda_1}),
\overrightarrow{\lambda_{\overline{1}}};
\lambda^{(4)},\cdots,\lambda^{(k)}).
\end{equation*}
For the second case, we can similarly construct $\widehat{O}_2$ which fixes all the parts
$\not\equiv -1 \pmod{3}$ and alternate the parts $\equiv -1 \pmod{3}$ in $\lambda^{(1)}$, $\lambda^{(2)}$, $\lambda^{(3)}$.
Let
\begin{equation}\label{HLW-O}
\widehat{O}=
\begin{cases}
\widehat{O}_1&\text{if $\lambda \in \mathcal{R}^{3n+j}$ fulfills the condition in the first case};\\
\widehat{O}_2&\text{if $\lambda \in \mathcal{R}^{3n+j}$ fulfills the conditions in the second case}.
\end{cases}
\end{equation}
Summing up the above suggests that the set $V_{t,k}^{3n+j}$
can be organized into orbits, with each orbit consists of 3 distinct members
\[
\Lambda(\lambda),\widehat{O}(\Lambda(\lambda)),\widehat{O}^2(\Lambda(\lambda)),
\]
and each of their $V_R$-rank
has distinct residue mod 3.
Clearly, the
total number of such orbits is $\frac{v_{t,k}(3n+j)}{3}$, and this
summarizes the combinatorial proof of $v_{t,k}(3n+j)\equiv 0 \pmod{3}$.
Since $\Lambda$ is a bijection and
by the definition of $V_R$-rank,
we conclude that the set $\mathcal{R}^{3n+j}$ can be organized into orbits with each orbit consists of three distinct members
\[
\lambda,\Lambda^{-1}\{\widehat{O}[\Lambda(\lambda)]\},
\Lambda^{-1}\{\widehat{O}^2[\Lambda(\lambda)]\},
\]
and each element of the orbit has a distinct $V_R$-rank $\pmod{3}$. This completes the proof.
\qed

%%%%%%%%%%%%%%%%%%%%%%%%%%%%%%%%%%%%%%%%%%%%%%%
\section{Proof of congruence \eqref{pd-3}}Let us give a quick overview of the notions for the 2-core and 2-quotient of an ordinary partition which are crucial in our proofs. For a lucid exposition see Schmidt \cite{Schmidt-2002}, and a full account of these topics can be found in \cite[pp. 75--85]{Kerber-1981}.

Firstly, recall that the Ferrers graph of the partition $\lambda$ is a set of coordinates in the bottom right quadrant of the plane where the $i$-th row contains $\lambda_i$ dots. We denote $\lambda'$ the conjugate of $\lambda$, which is the partition
whose graph is obtained by reflecting the Ferrers graph $\lambda$ about the main diagonal. For example,
we give $\lambda=4+4+2+2+1$ and its conjugate partition $\lambda'=5+4+2+2$ in Figure 1.
\vskip 3mm

\begin{figure}[htb]
\setlength{\unitlength}{1mm}
\begin{center}
\begin{picture}(41,30)(-0.5,-0.5)
\put(-24,0){$\bullet$}
\put(-24,6){$\bullet$}
\put(-24,12){$\bullet$}
\put(-24,18){$\bullet$}
\put(-24,24){$\bullet$}
%%%
\put(-18,6){$\bullet$}
\put(-18,12){$\bullet$}
\put(-18,18){$\bullet$}
\put(-18,24){$\bullet$}
%%%
\put(-12,18){$\bullet$}
\put(-12,24){$\bullet$}
%%%
\put(-6,18){$\bullet$}
\put(-6,24){$\bullet$}
%%%%%%%%%%%%%%%%%%%%%%%%%%%%%%%%%%%%%%%%%%%%%%%%%%
%%%%%%%%%%%%%%%%%%%%%%%%%%%%%%%%%%%%%%%%%
\put(36,6){$\bullet$}
\put(36,12){$\bullet$}
\put(36,18){$\bullet$}
\put(36,24){$\bullet$}
%%%
\put(42,6){$\bullet$}
\put(42,12){$\bullet$}
\put(42,18){$\bullet$}
\put(42,24){$\bullet$}
%%%
\put(48,18){$\bullet$}
\put(48,24){$\bullet$}
%%%
%%%
\put(54,18){$\bullet$}
\put(54,24){$\bullet$}
%%%
\put(60,24){$\bullet$}
\end{picture}
\end{center}
\caption{Ferrers graph of partitions
$\lambda=(4,4,2,2,1)$ and $\lambda'=(5,4,2,2)$.}
\label{fig:1}
\end{figure}
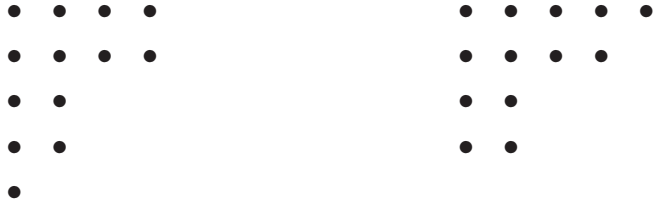
The 2-core of a partition $\lambda$ can be
formed as follows. Draw the Ferrers graph of $\lambda$ and successively remove 2-hooks $\bullet\ \  \bullet$ or  $\begin{matrix} \bullet \\ \bullet \end{matrix}$,
leaving a valid Ferrers graph at each stage, until
no more 2-hooks can be removed. It is easily to check the Ferrers graph that remains is
independent of the order in which the 2-hooks are removed and the remaining
graph is the 2-core of $\lambda$, denoted by $\lambda_{2c}$ in this paper.
Clearly, the 2-core of a partition is always a staircase partition, i.e., a partition of the form
\[
{n+1 \choose 2}=n+(n-1)+\cdots+1.
\]
The 2-quotient of a partition $\lambda$ is a bipartition $(\lambda_q^{(1)},\lambda_q^{(2)})\in \mathcal{P}_{-2}$ which can be constructed as follows.
Alternative write 0s and 1s on the coordinates in the Ferrers graph for $\lambda$, starting with a 0 in the top, left-most coordinate. Draw a horizontal line
through any row that ends in 0 and a vertical line through any column that
ends in 1. The Ferrers graph of $\lambda_q^{(1)}$ consists of the coordinates that are contained in both a vertical and a horizontal line; these coordinates are then pushed toward the northwest, in order to be justified with respect to the left and top, to form the Ferrers graph of $\lambda_q^{(1)}$. For $\lambda_q^{(2)}$, we carry out the same procedure with
the roles of rows and columns reversed.
An essential result given in \cite[3, p. 83, Theorem 2.7.30]{Nakayama-1940} points out that the triple
\[
(\lambda_{2c},\lambda_q^{(1)},\lambda_q^{(2)}),
\lambda_{2c} \in \mathcal{S}, \lambda_q^{(1)},\lambda_q^{(2)} \in \mathcal{P}_{1_{\{0\}}}
\]
uniquely determines $\lambda$, with
\[
|\lambda|=|\lambda_{2c}|+2|\lambda_q^{(1)}|+2|\lambda_q^{(2)}|.
\]
For simplicity, in this paper we denote the bijection
\begin{equation}
\Phi:\mathcal{P}_{1_{\{0\}}}^n \rightarrow (\lambda_{2c},\lambda^{(1)},\lambda^{(2)}),
\end{equation}
where
$\lambda_{2c} \in \mathcal{S}, \lambda^{(1)},\lambda^{(2)} \in \mathcal{P}_{2_{\{0\}}},
|\lambda_{2c}|+|\lambda^{(1)}|+|\lambda^{(2)}|=n$.
Take $\lambda=4+4+2+2+1$ as an example.
It is easily verified that
\[
\Phi(\lambda)=(\lambda_{2c},\lambda^{(1)},\lambda^{(2)})
=(1,2,6+4).
\]

To establish combinatorial proofs of congruence \eqref{pd-3},
we need to give a new representation of a partition with designated summands by a quintuple of partitions.
\begin{lem}\label{Theta}
There is a bijection $\Lambda_{pd}: \mathcal{PD}^n\rightarrow V_{2,5}^n$, where
\begin{align}\label{vtk}
\nonumber V_{2,5}=&\left\{(\lambda^{(1)},\lambda^{(2)},\lambda^{(3)},\lambda^{(4)},\lambda^{(5)})|
\lambda^{(1)},\lambda^{(2)},\lambda^{(3)} \in \mathcal{P}_{2_{\{0\}}}, \lambda^{(4)} \in \mathcal{S}, \lambda^{(5)} \in \mathcal{D}_{3_{\{0\}}}
\right\}.
\end{align}
\end{lem}
\proof
Let $\lambda\in \mathcal{PD}^n$.
Suppose that $d$ is a magnitude which appears in $\lambda$ and there are $m_d$ parts equal to $d$ among
which the $i_d$-th (from left to right) part is designated. There are two cases.
\begin{itemize}
\item[(1)] If $i_d=1$, then move all the parts equal to $d$ (including the designated part) in $\lambda$ to the partition $\alpha$;
\item[(2)] If $i_d \neq 1$, then move $i_d$ parts equal to $d$ in $\lambda$ to $\beta$ and $(m_d-i_d)$ parts equal to $d$ in $\lambda$ to $\alpha$.
\end{itemize}
Clearly, $\alpha$ is an ordinary partition and $\beta$ is a partition with each part occurs at least twice.
Write $\beta$ as in the form of $(1^{m_1}2^{m_2}\ldots)$
where $m_d$ is the multiplicity of $d$. Since $m_d\neq1$ for any $d$, there are two cases.
\begin{itemize}
\item[(1)] If $m_d$ is even, add $\frac{m_d}{2}$ parts $2d$ in $\lambda^{(3)}$;
\item[(2)] If $m_d$ is odd, add one part $3d$ in $\lambda^{(5)}$ and $\frac{m_d-3}{2}$ parts $2d$ in $\lambda^{(3)}$.
\end{itemize}
One can see that
$\lambda^{(3)} \in \mathcal{P}_{2_{\{0\}}}, \lambda^{(5)}\in \mathcal{D}_{3_{\{0\}}}$ and $|\beta|=|\lambda^{(3)}|+|\lambda^{(5)}|$. Denote this bijection by $\Psi$.
Recollect the bijection $\Phi$ from an ordinary partition to triples $(\lambda_{2c},\lambda^{(1)},\lambda^{(2)})$,
and let $\Phi(\alpha)=(\lambda^{(4)},\lambda^{(1)},\lambda^{(2)})$.
Now we can combine bijections $\Delta$, $\Psi$ and $\Phi$ into a single bijection $\Lambda_{pd}$ given as
\begin{equation*}
\Lambda_{pd}(\lambda)=
(\lambda^{(1)},\lambda^{(2)},\lambda^{(3)},\lambda^{(4)},\lambda^{(5)}).
\end{equation*}
For instance, let
\[
\lambda=20+20+20'+4+4'+4+4+2'+2
+1+1+1+1+1+1+1'+1,
\]
be a partition of 88 with designated summands.
Naturally,
\[
\Delta(\lambda)=(\alpha,\beta)=(4+4+2+2+1,20+20+20+4+4+1+1+1+1+1+1+1).
\]
Then we have
\[
\Phi(\alpha)=(\lambda^{(4)},\lambda^{(1)},\lambda^{(2)})=(1,2,6+4),
\]
and
\[
\Psi(\beta)=(\lambda^{(3)},\lambda^{(5)})=
(8+2+2,60+3).
\]
Eventually, we obtain
\begin{equation*}
\Lambda_{pd}(\lambda)=
(2,6+4,8+2+2,1,60+3).
\end{equation*}
One can check $\Lambda_{pd}$ is a bijection,
and the generating functions of $\lambda^{(4)},\lambda^{(5)}$
can be given as
\[
\sum_{n=0}^{\infty}q^{n+1 \choose 2},
(-q^3;q^3)_{\infty},
\]
respectively. \qed
\emph{Proof of congruence \eqref{pd-3}.}
Combining
\[
{n+1 \choose 2}\equiv 0 \ \  \textrm{or} \ \ 1 \pmod 3,
\]
and the fact $|\lambda^{(5)}|\equiv 0 \pmod{3}$,
we can deduce
\begin{equation*}
\sum_{i=4}^5|\lambda^{(i)}|\not\equiv 2 \pmod3.
\end{equation*}
Applying Theorem \ref{HLW-ZZ} leads to a combinatorial proof of congruence \eqref{pd-3}.
\qed
For example, in Table \ref{table-1} we list out the 15 partitions of $5$ with designated summands.
Apparently, there are $5$ orbits and each orbit contains 3 elements.

\begin{table}[htb]
\centering
% table caption is above the table
\caption{The case for $n=5$.}
\vskip 2mm
\label{table-1}       % Give a unique label
% For LaTeX tables use
\begin{tabular}{llcr}
\noalign{\smallskip}
$\lambda$ & $(\lambda^{(1)},\lambda^{(2)},\lambda^{(3)},\lambda^{(4)},\lambda^{(5)})$ &$r_V(\lambda)$
&$orbit$\\
\noalign{\smallskip}\hline\noalign{\smallskip}
$5'$           & $\left(4,\emptyset,\emptyset,1,\emptyset\right)$          &$1$      & $O_1$\\
$4'+1'$        &  $\left(\emptyset,2,\emptyset,2+1,\emptyset\right)$       & $-1$    & $O_2$\\
$3'+2'$        &  $\left(2+2,\emptyset,\emptyset,1,\emptyset\right)$       & $2$     & $O_3$\\
$3'+1'+1$      &  $(2,2,\emptyset,1,\emptyset)$     & $0$    & $O_4$\\
$3'+1+1'$      &  $(2,\emptyset,2,1,\emptyset)$               & $1$     & $O_4$\\
$2'+2+1'$      &  $(\emptyset,4,\emptyset,1,\emptyset)$     & $-1$     & $O_1$\\
$2+2'+1'$      &  $(\emptyset,\emptyset,4,1,\emptyset)$               & $0$     & $O_1$\\
$2'+1'+1+1$    &  $(2,\emptyset,\emptyset,2+1,\emptyset)$   & $1$    & $O_2$\\
$2'+1+1'+1$     &  $(\emptyset,\emptyset,2,2+1,\emptyset)$            & $0$     & $O_2$\\
$2'+1+1+1'$     &  $(\emptyset,2,\emptyset,\emptyset,3)$            & $-1$     & $O_5$\\
$1'+1+1+1+1$     &  $(\emptyset,2+2,\emptyset,1,\emptyset)$       & $-2$  & $O_3$\\
$1+1'+1+1+1$     &  $(\emptyset,2,2,1,\emptyset)$                 & $-1$  & $O_4$\\
$1+1+1'+1+1$     &  $(2,\emptyset,\emptyset,\emptyset,3)$                 & $1$  & $O_5$\\
$1+1+1+1'+1$     &  $(\emptyset,\emptyset,2+2,1,\emptyset)$                 & $0$   & $O_3$\\
$1+1+1+1+1'$     &  $(\emptyset,\emptyset,2,\emptyset,3)$       & $0$  & $O_5$\\
\noalign{\smallskip}\hline
\end{tabular}
\end{table}

%%%%%%%%%%%%%%%%%%%%%%%%%%%%%%%%%%%%%%%%%%
%%%%%%%%%%%%%%%%%%%%%%%%%%%%%%%%%%%%%%%%%%%%

\section{Proof of congruence \eqref{a-3}}

In \cite{Kim-2010}, Kim interpreted $a(n)$ as the number of 2-color restricted partitions of $n$ with colors $r$
and $b$ such that the color $b$ appears only in even parts.
For instance, there are 3 such partitions of 2:
\[
2^r,2^b,1^r+1^r.
\]
Invoking Theorem \ref{HLW-ZZ}, the combinatorial proof of congruence \eqref{a-3} follows immediately according to the lemma given below.
\begin{lem}\label{Lambda}
There is a bijection $\Lambda_a: \mathcal{A}^n\rightarrow V_{2,4}^n$, where
\begin{align}\label{vtk-cubic}
\nonumber V_{2,4}=&\left\{(\lambda^{(1)},\lambda^{(2)},\lambda^{(3)},\lambda^{(4)})|
\lambda^{(1)},\lambda^{(2)},\lambda^{(3)} \in \mathcal{P}_{2_{\{0\}}}, \lambda^{(4)} \in \mathcal{S}
\right\}.
\end{align}
\end{lem}
\proof Let $\lambda \in \mathcal{A}^n$. Clearly, we can split $\lambda$ into the bipartition $(\alpha,\lambda^{(3)})$ in accordance with the color of the parts, namely, where $\alpha$
is consisted of the parts colored $r$ and $\lambda^{(3)}$ is composed of the parts colored $b$. It is obvious that
$\alpha \in \mathcal{P}_{1_{\{0\}}}$ and $\lambda^{(3)} \in \mathcal{P}_{2_{\{0\}}}$.
Recollect the bijection $\Phi$ from ordinary partition to triples $\left(\lambda_{2c},\lambda_0,\lambda_1\right)$, and
let $\Phi(\alpha)=(\lambda^{(4)},\lambda^{(1)},\lambda^{(2)})$.
Summing up the above, we reach Lemma \ref{Lambda}, as claimed. \qed

The proof of congruence \eqref{a-3} is similar to that of congruence \eqref{pd-3}, and hence it is omitted. We just give an example in Table \ref{table-2}.

\begin{table}[htb]
\centering
% table caption is above the table
\caption{The case for $n=5$.}
\vskip 2mm
\label{table-2}       % Give a unique label
% For LaTeX tables use
\begin{tabular}{llcr}
\noalign{\smallskip}
$\lambda$ & $\left(\lambda^{(1)},\lambda^{(2)},\lambda^{(3)},\lambda^{(4)}\right)$ & $r_V(\lambda)$
& $orbit$\\
\noalign{\smallskip}\hline\noalign{\smallskip}
$5^r$           & $\left(4,\emptyset,\emptyset,1\right)$          &$1$      & $O_1$\\
$4^r+1^r$        &  $\left(\emptyset,2,\emptyset,2+1\right)$       & $-1$    & $O_2$\\
$4^b+1^r$        &  $\left(\emptyset,\emptyset,4,1\right)$       & $0$     & $O_1$\\
$3^r+2^r$      &  $(2+2,\emptyset,\emptyset,1)$     & $2$    & $O_3$\\
$3^r+2^b$      &  $(2,\emptyset,2,1)$               & $1$     & $O_4$\\
$3^r+1^r+1^r$      &  $(2,2,\emptyset,1)$     & $0$     & $O_4$\\
$2^r+2^r+1^r$      &  $(\emptyset,4,\emptyset,1)$               & $-1$     & $O_1$\\
$2^b+2^r+1^r$    &  $(\emptyset,\emptyset,2,2+1)$   & $0$    & $O_2$\\
$2^b+2^b+1^r$     &  $(\emptyset,\emptyset,2+2,1)$            & $0$     & $O_3$\\
$2^r+1^r+1^r+1^r$     &  $(2,\emptyset,\emptyset,2+1)$            & $1$     & $O_2$\\
$2^b+1^r+1^r+1^r$     &  $(\emptyset,2,2,1)$       & $-1$  & $O_4$\\
$1^r+1^r+1^r+1^r+1^r$     &  $(\emptyset,2+2,\emptyset,1)$                 & $-2$  & $O_3$\\
\noalign{\smallskip}\hline
\end{tabular}
\end{table}
\section{Proofs of congruences \eqref{pod-3} and \eqref{o-3}}

We first prove congruence \eqref{pod-3}.
\begin{lem}\label{Omega}
There is a bijection $\Lambda_{pod}: \mathcal{POD}_{-2}^n\rightarrow V_{2,4}^n$, where
\begin{align}
\nonumber V_{2,4}=&\left\{(\lambda^{(1)},\lambda^{(2)},\lambda^{(3)},\lambda^{(4)})|
\lambda^{(1)},\lambda^{(2)},\lambda^{(3)} \in \mathcal{P}_{2_{\{0\}}}, \lambda^{(4)} \in \mathcal{S}_{t-odd}
\right\}.
\end{align}
\end{lem}
\proof
To avoid conflicts on the symbols, in the proof of Lemma \ref{Omega}, we write a partition as $\lambda=(\lambda_1,\ldots, \lambda_l)$. For a given bipartition $(\lambda_{pod}^{(1)},\lambda_{pod}^{(2)}) \in \mathcal{POD}_{-2}^n$, take out all the even parts in $\lambda_{pod}^{(1)}$
to form $\lambda^{(1)}$, and all the even parts in $\lambda_{pod}^{(2)}$
to form $\lambda^{(2)}$.
Clearly, the remaining parts in
$(\lambda_{pod}^{(1)},\lambda_{pod}^{(2)})$ form a bipartition
$(\mu^{(1)},\mu^{(2)})$ with
$\mu^{(1)}, \mu^{(2)} \in \mathcal{D}_{2_{\{1\}}}$.
Hence it suffices to show that
there is a bijection between bipartitions $(\mu^{(1)},\mu^{(2)})$ of $n$ and bipartitions $(\lambda^{(3)},\lambda^{(4)})$ of $n$ where $\lambda^{(3)} \in \mathcal{P}_{2_{\{0\}}}$ and $\lambda^{(4)} \in \mathcal{S}_{t-odd}$.
A graceful modified version of the Wright map $\varphi$ constructed by Seo and Yee \cite{Seo-Yee-2017} is significant for our proof in this part, and we denote this modified version by $\varphi_m$ (see \cite{Wright-1965,Yee-2015} for details of the Wright map).
Recall that a Frobenius symbol of $n$ is a two-rowed array\cite{Yee-2003}
\[F=\left(
\begin{array}{ccccc}
  a_1 & a_2 & \cdots & a_l \\
b_1 & b_2 & \cdots & b_l
\end{array}\right),
\]
where $a_1>a_2>\ldots>a_l\geq 0$, $b_1>b_2>\ldots>b_l\geq 0$ and $n=\sum_{i=1}^l (a_i+b_i)+l$.
Given the Ferrers graph of an ordinary partition,
$a_i$ form rows to the right of the diagonal, and $b_i$ form columns
below the diagonal. As a result, there is a natural  bijection between
the Frobenius symbols of $n$ and the ordinary partitions of $n$.
We are now ready to apply the bijection $\varphi_m$ for
$(\mu^{(1)}, \mu^{(2)})$, namely
\begin{align*}
\mu^{(1)}&=(2a_1+1,2a_2+1,\ldots,2a_{l+m}+1),\\[5pt]
\mu^{(2)}&=(2b_1+1,2b_2+1,\ldots,2b_{l}+1),
\end{align*}
where $a_1>a_2>\cdots>a_{l+m}\geq0$ and $b_1>b_2>\cdots>b_l\geq 0$.
Suppose that $\mu$ and $\nu$ are two partitions. For convenience, let $\mu\cup\nu$ be the partition consisting of all the parts of $\mu$ and $\nu$, and denote $2\mu$ as the partition
whose parts are $2$ times each part of $\mu$. For
instance, let $\mu=(9,6,6,2,1)$. Accordingly, $2\mu=(18,12,12,4,2)$.
For $m\geq 0$, we construct a Forbenius symbol
\begin{align*}
\mu=\left(
\begin{array}{llll}
a_{1+m} & a_{2+m} &\cdots  &a_{l+m} \\[5pt]
 b_1 & b_2 & \cdots & b_l
\end{array}\right)
\end{align*}
and a partition $\nu=(a_1-m+1,a_2-m+2,\ldots,a_m)$.
We define $\varphi_m(\mu^{(1)},\mu^{(2)})=(\pi,\triangle)$, where $\pi=2(\mu\cup\nu)$ and $\triangle=(2(m-1)+1,2(m-2)+1,\ldots,3,1)$. When $m=0$, we see that $\triangle=\emptyset$.
For $m<0$, a Forbenius symbol
\begin{eqnarray*}
\mu=\left(
\begin{array}{llll}
b_{1-m} & b_{2-m} & \cdots & b_l \\[5pt]
a_{1} & a_{2} &\cdots  &a_{l+m}
\end{array}\right)
\end{eqnarray*}
and a partition $\nu=(b_1+m+1,b_2+m+2,\ldots,b_{-m})$ are built.
We define $\varphi_m(\mu^1,\mu^2)=(\pi,\triangle)$, where $\pi=2(\mu\cup\nu)'$ and $\triangle=(2(-m-1)+1,2(-m-2)+1,\ldots,3,\overline{1})$.
For instance, give $(\mu^{(1)},\mu^{(2)})$ as
\begin{align*}
\mu^{(1)}&=(9,7,3)=(2\times4+1,2\times3+1,2\times1+1),\\[5pt]
\mu^{(2)}&=(17,15,11,7,3,1)=(2\times8+1,2\times7+1,2\times5+1,2\times3+1,2\times1+1,2\times0+1).
\end{align*}
In light of the modified Wright's map $\varphi_m$, we get
\begin{align*}
&\mu=
\left(
\begin{array}{lll}
b_{4} & b_{5} & b_{6} \\[5pt]
a_{1} & a_{2} & a_{3}
\end{array}\right)=
\begin{pmatrix}
3 & 1 & 0\\ 4 & 3 & 1
\end{pmatrix}
=(4,3,3,3,2),\\[3pt]
&\nu=(8-3+1,7-3+2,5-3+3)=(6,6,5).
\end{align*}
Hence we obtain
\begin{align*}
&\varphi_m(\mu^{(1)},\mu^{(2)})=(\pi,\triangle)=(2(\mu\cup \nu)',((3-1)\times 2+1,(2-1)\times 2+1,\overline{1}))\\[3pt]
=&(2(6,6,5,4,3,3,3,2)',(5,3,\overline{1}))
=((16,16,14,8,6,4),(5,3,\overline{1})).
\end{align*}
It can be checked that $|\mu^1|+|\mu^2|=|\pi|+|\triangle|,\pi \in \mathcal{P}_{2_{\{0\}}}, \triangle \in \mathcal{S}_{t-odd}$ and the generating function of $\triangle$ can be given as
\[
q^0+2\sum_{n=1}^{\infty}q^{n^2} =\sum_{n=-\infty}^{\infty}q^{n^2}.
\]
By now we complete the proof of Lemma \ref{Omega}.
\qed
\emph{Proof of congruence \eqref{pod-3}.}
Since
$n^2\not\equiv 2 \pmod 3$, invoking Theorem \ref{HLW-ZZ}, we give a combinatorial proof of congruence \eqref{pod-3}. Here we give an example in Table \ref{table-3}.\qed

\begin{table}[htb]
\centering
% table caption is above the table
\caption{The case for $n=5$.}
\vskip 2mm
\label{table-3}       % Give a unique label
% For LaTeX tables use
\begin{tabular}{llcr}
\noalign{\smallskip}
$\lambda$ & $\left(\lambda^{(1)},\lambda^{(2)},\lambda^{(3)},\lambda^{(4)}\right)$  & $r_V(\lambda)$
& $orbit$\\
\noalign{\smallskip}\hline\noalign{\smallskip}
$(5,\emptyset)$
&$\left(\emptyset,\emptyset,4,1\right)$
&$0$
&$O_1$\\
%%%%%%%%%%%%%%%%%%%%%%%%%%%%%%%%
$(\emptyset,5)$
&$\left(\emptyset,\emptyset,2+2,\overline{1}\right)$
&$0$
&$O_2$\\
%%%%%%%%%%%%%%%%%%%%%%%%%%%%%%%%
$(4+1,\emptyset)$
&$\left(4,\emptyset,\emptyset,1\right)$
&$1$
&$O_1$\\
%%%%%%%%%%%%%%%%%%%%%%%%%%%%%%%%
$(\emptyset,4+1)$
&$\left(\emptyset,4,\emptyset,\overline{1}\right)$
&$-1$
&$O_3$\\
%%%%%%%%%%%%%%%%%%%%%%%%%%%%%%%%
$(4,1)$
&$\left(4,\emptyset,\emptyset,\overline{1}\right)$
&$1$
&$O_3$\\
%%%%%%%%%%%%%%%%%%%%%%%%%%%%%%%%
$(1,4)$
&$\left(\emptyset,4,\emptyset,1\right)$
&$-1$
&$O_1$\\
%%%%%%%%%%%%%%%%%%%%%%%%%%%%%%%%
$(3+2,\emptyset)$
&$\left(2,\emptyset,2,1\right)$
&$1$
&$O_4$\\
%%%%%%%%%%%%%%%%%%%%%%%%%%%%%%%%
$(\emptyset,3+2)$
&$\left(\emptyset,2,2,\overline{1}\right)$
&$-1$
&$O_5$\\
%%%%%%%%%%%%%%%%%%%%%%%%%%%%%%%%
$(3,2)$
&$\left(\emptyset,2,2,1\right)$
&$-1$
&$O_4$\\
%%%%%%%%%%%%%%%%%%%%%%%%%%%%%%%%
$(2,3)$
&$\left(2,\emptyset,2,\overline{1}\right)$
&$1$
&$O_5$\\
%%%%%%%%%%%%%%%%%%%%%%%%%%%%%%%%
%%%%%%%%%%%%%%%%%%%%%%%%%%%%%%%%
$(3+1,1)$
&$\left(\emptyset,\emptyset,2+2,1\right)$
&$0$
&$O_6$\\
%%%%%%%%%%%%%%%%%%%%%%%%%%%%%%%%
$(1,3+1)$
&$\left(\emptyset,\emptyset,4,\overline{1}\right)$
&$0$
&$O_3$\\
%%%%%%%%%%%%%%%%%%%%%%%%%%%%%%%%
$(2+2+1,\emptyset)$
&$\left(2+2,\emptyset,\emptyset,1\right)$
&$2$
&$O_6$\\
%%%%%%%%%%%%%%%%%%%%%%%%%%%%%%%%
$(\emptyset,2+2+1)$
&$\left(\emptyset,2+2,\emptyset,\overline{1}\right)$
&$-2$
&$O_2$\\
%%%%%%%%%%%%%%%%%%%%%%%%%%%%%%%%
$(2+2,1)$
&$\left(2+2,\emptyset,\emptyset,\overline{1}\right)$
&$2$
&$O_2$\\
%%%%%%%%%%%%%%%%%%%%%%%%%%%%%%%%
$(1,2+2)$
&$\left(\emptyset,2+2,\emptyset,1\right)$
&$-2$
&$O_6$\\
%%%%%%%%%%%%%%%%%%%%%%%%%%%%%%%%
$(2+1,2)$
&$\left(2,2,\emptyset,1\right)$
&$0$
&$O_4$\\
%%%%%%%%%%%%%%%%%%%%%%%%%%%%%%%%
$(2,2+1)$
&$\left(2,2,\emptyset,\overline{1}\right)$
&$0$
&$O_5$\\
%%%%%%%%%%%%%%%%%%%%%%%%%%%%%%%%
\noalign{\smallskip}\hline
\end{tabular}
\end{table}
Congruence \eqref{o-3} can be justified in the same manner according to the Wright bijection, and hence it is omited.
\begin{rem}
Let $O$ be an injective
map on $V_{t,k}$ defined as
\begin{align*}
&O((\lambda^{(1)},\lambda^{(2)},\lambda^{(3)},\lambda^{(4)}\cdots,\lambda^{(k)}))
=(\lambda^{(3)},\lambda^{(1)},\lambda^{(2)},\lambda^{(4)}\cdots,\lambda^{(k)}).
\end{align*}
It is obvious that the fixed-point set of the injective map $O$ is the subset of
$V_{t,k}$ subject to the restriction that $\lambda^{(1)}=\lambda^{(2)}=\lambda^{(3)}$. In light of this observation, if $\mathcal{R}$ fulfills the conditions of Theorem \ref{HLW-ZZ}, we can give another method to organize the set $\mathcal{R}^{3n+j}$ into orbits under the injective map $O$, such that each orbit consists of 3 distinct members.
\end{rem}
\begin{rem}
Denote
\begin{align*}
\mathcal{R}_g=&\{\lambda| \lambda\  \textrm{is a partition subjects to given restrictions}\}.
\end{align*}
Let
\begin{align}\label{vRk}
\nonumber V_{R_g,k}=&\left\{(\lambda^{(1)},\lambda^{(2)},\lambda^{(3)},\cdots,\lambda^{(k)})|
\lambda^{(1)},\lambda^{(2)},\lambda^{(3)} \in \mathcal{R}_g, \ \textrm{for}\  i>3, \right. \\[3pt]
&\left.\lambda^{(i)} \ \textrm{is an ordinary partition or a restricted partition}
\right\}
\end{align}
subjecting to the restriction that
$\widehat{O}$ is an injective
map on $V_{R_g,k}$,
where $\widehat{O}$ is defined by \eqref{HLW-O}.
Clearly, if we replace $V_{t,k}^n$ by $V_{R_g,k}^n$ in Theorems \ref{HLW-ZZ}, the theorem still holds. Based on this observation one can easily give combinatorial proofs of the congruences mod 3 given by Toh \cite{Toh-2012}, and we omit the details.
\end{rem}

\noindent{\bf Acknowledgments.}

The author was supported by the Scientific Research Foundation of Nanjing Institute of Technology.

%%%%%%%%%%%%%%%%%%%%%%%%%%%%%%%%%%%%%%%%%%%%%%%%%%%%%%%%%%%%%%%%


\begin{thebibliography}{}

\bibitem{Andrews-Garvan-1988}Andrews, G. E., Garvan, F. G.:
Dyson's crank of a partition,  Bull. Amer. Math. Soc.
{\bf 18}, 167--171 (1988)


\bibitem{Andrews-Lewis-Lovejoy-2002}
Andrews, G. E., Lewis, R. P.,  Lovejoy, J.: Partitions with designated summands, Acta Arith. {\bf 105}, 51--66 (2002)

\bibitem{Chan-2010} Chan, H.-C.: Ramanujan's cubic continued fraction and an analog of his ``most beautiful identity'', Int. J. Number
Theory {\bf 06}, 673--680 (2010)

\bibitem{Chen-Ji-Jin-Shen-2013}Chen, W. Y. C., Ji, K. Q., Jin, H-T.,
Shen, E. Y. Y.: On the number of partitions with designated summands, J. Number Theory {\bf 133}, 2929--2938 (2013)

\bibitem{Chen-Lin-2011} Chen, W. Y. C., Lin, B. L. S.:  Congruences for bipartitions with odd parts distinct,
Ramanujan J. {\bf 25}, 277--293 (2011)

\bibitem{Chen-Lin-2012} Chen, W. Y. C., Lin, B. L. S.:  Arithmetic properties of overpartition pairs,
Acta Arithmetica (2012)

\bibitem{Fu-Tang-2017} Fu, S. S., Tang D. Z.: Multiranks and classical Theta functions, Int. J. Number
Theory {\bf 14}, 549--566 (2017).

\bibitem{Garvan-1988} Garvan, F. G.: New combinatorial interpretations of Ramanujan's partition congruences mod $5$, $7$, $11$, Trans. Amer. Math. Soc. {\bf 305}, 47--77 (1988)

\bibitem{Garvan-Kim-Stanton-1990} Garvan, F. G., Kim, D., Stanton, D.: Cranks and $t$-cores, Invent. Math. {\bf 101}, 1--17 (1990)

\bibitem{Hirschhorn-Sellers-2005} Hirschhorn, M. D., Sellers, J. A.: Arithmetic relations for overpartitions, J. Comb. Math. Comb. Comp. {\bf 53},  65--73 (2005)

\bibitem{Kerber-1981} James, G., Kerber, A.: The Representation Theory of the Symmetric Group, Reading, MA: Addison-Wesley, 1981.

\bibitem{Kim-2010} Kim, B.: An analog of crank for a certain kind of partition function arising from the cubic continued fraction, Acta Arithmetica {\bf 148}, 1--19 (2010)

%\bibitem{Lin-2018} Lin, B. L. S.:  The number of tagged parts over the
%partitions with designated summands, J. Number Theory, {\bf 184} 216--234 (2018)

\bibitem{Nakayama-1940} Nakayama, T.: On some modular properties of irreducible representations of a symmetric group, I, II, Japan. J. Math. {\bf 17}, 165--184, 411--423 (1940)

\bibitem{Ramanujan-1919} Ramanujan, S.: Some properties of $p(n)$, the number of partitons of $n$, Proc. Cambridge
Philos. Soc. {\bf 19}, 214--216 (1919)

\bibitem{Schmidt-2002} Schmidt, F.: Integer partitions and binary trees, Advances in Applied Mathematics, (2002).

\bibitem{Seo-Yee-2017}Seo, S., Yee, A. J.:, Overpartitions and singular overpartitions, analytic number theory, modular forms and $q$-Hypergeometric series: In Honor of Krishna Alladi's 60th Birthday, University of Florida, Gainesville, 693--711 (2016)

\bibitem{Toh-2012} Toh, P. C.: Ramanujan type identities and congruences for partition pairs, Discrete Math. {\bf 312}, 1244--1250 (2012)

\bibitem{Wang-2017} Wang, L. Q.: Congruences modulo powers of 11 for some partition functions, Proceedings of the American Mathematical Society {\bf 146}, 1515--1528 (2017)

\bibitem{Wang-Yee-2020} Wang, C., Yee, A. J.: Truncated Hecke-Rogers type series, Adv. Math. {\bf 365}, 51--70 (2020)

\bibitem{Wright-1965} Wright, E. M.: An enumerative proof of an identity of Jacobi, J. London Math. Soc. {\bf 40}, 55--57 (1965)

\bibitem{Yee-2003} Yee, A. J.: Combinatorial proofs of generating function identities for F-partitions, J. Combin. Ser. A {\bf 102}, 217--228 (2003)

\bibitem{Yee-2015} Yee, A. J.: Truncated Jacobi triple product theorem, J. Combin. Ser. A {\bf 130}, 1--14 (2015)





\end{thebibliography}
\end{document}